\documentclass[reqno]{amsart}
\usepackage{amsmath,amssymb,amsthm,ascmac}
\usepackage[mathscr,mathcal]{eucal}
\usepackage{latexsym}
\usepackage{enumerate}
\usepackage[dvipdfmx]{graphicx}

\begin{document}

\allowdisplaybreaks

\theoremstyle{definition}
\newtheorem{defi}{\textbf{Definition}}[section]
\newtheorem{thm}[defi]{\textbf{Theorem}}
\newtheorem{lem}[defi]{\textbf{Lemma}}
\newtheorem{prop}[defi]{\textbf{Proposition}}
\newtheorem{cor}[defi]{\textbf{Corollary}}
\newtheorem{ex}[defi]{\textbf{Example}}
\newtheorem{rem}[defi]{\textbf{Remark}}
\newtheorem*{corr}{\textbf{Corollary}}

\theoremstyle{plain}
\newtheorem{maintheorem}{Theorem}
\newtheorem{theorem}{Theorem }[section]
\newtheorem{proposition}[theorem]{Proposition}
\newtheorem{mainproposition}{Proposition}
\newtheorem{lemma}[theorem]{Lemma}
\newtheorem{corollary}[theorem]{Corollary}
\newtheorem{maincorollary}{Corollary}
\newtheorem{claim}{Claim}
\renewcommand{\themaintheorem}{\Alph{maintheorem}}
\theoremstyle{definition} \theoremstyle{remark}
\newtheorem{remark}[theorem]{Remark}
\newtheorem{example}[theorem]{Example}
\newtheorem{definition}[theorem]{Definition}
\newtheorem{problem}{Problem}
\newtheorem{question}{Question}
\newtheorem{exercise}{Exercise}

\newtheorem*{subject}{�ړI}
\newtheorem*{mondai}{Problem}
\newtheorem{lastpf}{�ؖ�}
\newtheorem{lastpf1}{proof of theorem2.4}
\renewcommand{\thelastpf}{}

\newtheorem{last}{Theorem}

\renewcommand{\thelast}{}
\renewcommand{\proofname}{\textup{Proof.}}

\renewcommand{\theequation}{\arabic{section}.\arabic{equation}}
\makeatletter
\@addtoreset{equation}{section}

\title[LDP for linear mod 1 transformations]
{Large deviation principle for linear mod 1 transformations}
\author[Y. M. Chung]{Yong Moo Chung}
\address{Department of Applied Mathematics \\
Hiroshima University \\
Higashi-Hiroshima 739-8527, JAPAN
}
\email{chung@amath.hiroshima-u.ac.jp}
\author[K. Yamamoto]{Kenichiro Yamamoto}
\address{Department of General Education \\
Nagaoka University of Technology \\
Nagaoka 940-2188, JAPAN
}
\email{k\_yamamoto@vos.nagaokaut.ac.jp}


\date{}

\maketitle
\large

\begin{abstract}
For $0\le\alpha<1$ and $\beta>2$, we consider a linear mod 1 transformation on a unit interval;
$x\mapsto \beta x+\alpha$ (${\rm mod }\ 1$), and prove that it satisfies the level-2 large deviation
principle with the unique measure of maximal entropy.
For the proof, we use the density of periodic measures and Hofbauer's Markov Diagram.
\end{abstract}

\section{Introduction}

Let $X$ be a metrizable space and $T\colon X\to X$ be a Borel measurable map.
We denote by $\mathcal{M}(X)$ the set of all Borel probability measures on $X$ endowed with
the weak${}^{\ast}$-topology,
by $\mathcal{M}_T(X)\subset \mathcal{M}(X)$ the set of all $T$-invariant ones and by
$\mathcal{M}^e_T(X)\subset\mathcal{M}_T(X)$ the set of ergodic ones.
We say that $(X,T)$ satisfies a \textit{(level-2) large deviation principle}
with a \textit{reference measure} $m\in\mathcal{M}(X)$
if there exists a lower semicontinuous function $\mathcal{J}\colon \mathcal{M}(X)\to [0,\infty]$, called
a \textit{rate function} such that
$$\limsup_{n\rightarrow\infty}\frac{1}{n}\log m(\{x\in X:\delta_n^T(x)\in\mathcal{K}\})\le
-\inf_{\mathcal{K}}\mathcal{J}$$
holds for any closed set $\mathcal{K}\subset \mathcal{M}(X)$ and
$$\liminf_{n\rightarrow\infty}\frac{1}{n}\log m(\{x\in X:\delta_n^T(x)\in\mathcal{U}\})\ge
-\inf_{\mathcal{U}}\mathcal{J}$$
holds for any open set $\mathcal{U}\subset\mathcal{M}(X)$.
Here $\delta_n^T\colon X\to\mathcal{M}(X)$ is defined by $\delta_n^T(x):=\frac{1}{n}\sum_{j=0}^{n-1}\delta_{T^j(x)}$
and $\delta_y$ stands for the Dirac mass at the point $y\in X$.
We say that $\mu\in\mathcal{M}(X)$ is a
\textit{periodic measure} if there exist $x\in X$ and $n>0$ such that
$T^n(x)=x$ and $\mu=\delta_n^T(x)$ hold. Then it is clear that $\mu\in\mathcal{M}_T^e(X)$.
We denote by $\mathcal{M}_T^p(X)\subset\mathcal{M}_T^e(X)$
the set of all periodic measures on $X$.

Hereafter, let $X=[0,1]$ be a unit interval and $T\colon X\to X$ be a
\textit{piecewise monotonic map}, i.e.,
there exist integer $k>1$ and $0=a_0<a_1<\cdots<a_k=1$ such that
$T|_{(a_{j-1},a_j)}$ is continuous and strictly monotone for
each $1\le j\le k$.
Throughout of this paper, we further assume the following conditions for a piecewise monotonic map $T$.
\begin{itemize}
\item
$T$ is \textit{transitive}, i.e., there exists a point $x\in [0,1]$ whose forward orbit
$\{T^n(x):n\ge 0\}$ is dense in $[0,1]$;

\item
The topological entropy $h_{\rm top}(T)$ of $T$ is positive.
\end{itemize}
Then it is proved in \cite{H2} that there exists a unique measure of maximal entropy.
In this paper, we investigate large deviations for piecewise monotonic maps
with the unique measure of maximal entropy as reference.
In such situation, it was shown in \cite{T1,Yo} that the large deviation principle holds for
maps with the specification property (including Markov maps and continuous maps).
After that, Pfister and Sullivan extended the results to more general maps beyond specification
including all $\beta$-transformations by using certain specification-like property, called the
approximate product property (\cite{PS}).

The aim of this paper is to extend large deviations results to linear mod 1 transformations.
A linear mod 1 transformation $T_{\alpha,\beta}\colon X\to X$ with $\beta>1$ and
$0\le \alpha<1$ was introduced by Parry (\cite{P}) and defined by
\begin{equation}
\label{d-mod-1}
T_{\alpha,\beta}(x)=
\left\{
\begin{array}{ll}
\beta x+\alpha\ (\text{mod}\ 1) & (x\in [0,1)),\\
\lim_{x\rightarrow 1-0}T(x) & (x=1).
\end{array}
\right.
\end{equation}
The map $T_{\alpha,\beta}$ is said to be a $\beta$-transformation if $\alpha=0$.
Throughout of this paper, we focus on linear mod 1 transformations with $\beta>2$.
Then it is straightforward to see that $T_{\alpha,\beta}$ is transitive and
hence there is a unique measure of maximal entropy.
On the other hand,
for Lebesgue almost every parameter $(\alpha,\beta)\in [0,1)\times (2,\infty)$,
the map $T_{\alpha,\beta}$ does not satisfy the specification property (\cite{Bu2}).
Moreover, it is unknown any suitable specification-like
condition to show the large deviation principle
for linear mod 1 transformations.
As far as we know, our main result is the first large deviations result
for such maps.

\begin{maintheorem}
\label{mod1}
Fix $0\le \alpha<1$ and $\beta>2$. Let $T=T_{\alpha,\beta}\colon X\to X$ be as above
and $m$ be the unique measure of maximal entropy.
Then $(X,T)$ satisfies the large deviation principle with $m$ and the rate function
$\mathcal{J}\colon\mathcal{M}(X)\to [0,\infty]$ is given by
$$\mathcal{J}(\mu)=
\left\{
\begin{array}{ll}
h_{\rm top}(T)-h_{\mu}(T) & (\mu\in \mathcal{M}_T(X)),\\
\infty & (\text{otherwise}).
\end{array}
\right.
$$
Here $h_{\mu}(T)$ denotes the metric entropy of $\mu\in\mathcal{M}_T(X)$.
\end{maintheorem}

For a piecewise monotonic map $T$, we consider the following condition
which is strictly stronger than the transitivity, but still holds for all linear mod 1 transformations with $\beta>2$:
\begin{equation}
\label{irreducible}
\begin{split}
& \text{For any open interval }I\subset [0,1],\text{ there exist an positive integer }\tau,\\
& \text{integers }1\le i_0,\ldots,i_{\tau-1}\le k\text{ and an open subinterval }L\subset I\text{ so that }\\
& T^j(L)\subset (a_{i_j-1},a_{i_j})\text{ for all }0\le j\le \tau-1\text{ and }T^{\tau}(L)=(0,1).
\end{split}
\end{equation}

Theorem \ref{mod1} is derived from
the following general result for piecewise monotonic maps.

\begin{maintheorem}
\label{interval}
Let $T\colon X\to X$ be a piecewise monotonic map satisfying (\ref{irreducible}) and
$m$ be the unique measure of maximal entropy.
Assume that $\mathcal{M}^p_T(X)$ is dense in $\mathcal{M}^e_T(X)$.
Then $(X,T)$ satisfies the large deviation principle with $m$ and the rate function
$\mathcal{J}\colon\mathcal{M}(X)\to [0,\infty]$ is given by
\begin{equation}
\label{rate2}
\mathcal{J}(\mu)=
\left\{
\begin{array}{ll}
h_{\rm top}(T)-h_{\mu}(T) & (\mu\in {\rm cl}(\mathcal{M}_T(X_T))), \\
\infty & (\text{otherwise}).
\end{array}
\right.
\end{equation}
Here $X_T:=\bigcap_{n=0}^{\infty}T^{-n}(\bigcup_{j=1}^k (a_{j-1},a_j))$
and ${\rm cl}(A)$ denotes the closure of a set $A$. 
\end{maintheorem}



\begin{rem}
In the equation (\ref{rate2}), ${\rm cl(\mathcal{M}_T(X_T))}$
cannot be replaced by $\mathcal{M}_T(X)$ in general.
Indeed, in \S4, we construct an example of piecewise monotonic maps $T$
to which Theorem \ref{interval} is applicable and $\mathcal{M}_T(X)$ does not coincide
with ${\rm cl}(\mathcal{M}_T(X_T))$ (see Example \ref{counter}).
\end{rem}

Theorem \ref{interval} has applications to the following important classes of piecewise monotonic maps:
\begin{itemize}
\item
All $\beta$-transformations (see \cite{G,Si} for instance).

\item
Linear mod 1 transformations with $\beta>2$ (see the proof of Theorem \ref{mod1} in \S3).

\item
$(-\beta)$-transformations with $(1+\sqrt{5})/2<\beta<2$ (see Example \ref{counter} for the definition and details).
\end{itemize}
In particular, we can and do apply Theorem \ref{interval} to linear mod 1 transformations
and $(-\beta)$-transformations.
To the best of our knowledge, this is the first result for large deviations of these transformations.





This paper is organized as follows:
In \S2, we establish our definitions and prepare several lemmas.
We show Theorems \ref{mod1} and \ref{interval} in \S3.
In \S4, we give other applications of Theorem \ref{interval}.

\section{Preliminaries}
\subsection{Symbolic dynamics}
Let $A$ be a countable set.
Denote by $A^{\mathbb N}$ the one-sided  infinite product of $A$ 
equipped with the product topology of the discrete topology of $A$.
Let $\sigma$ be the shift map on $A^{\mathbb N}$ 
 (i.e.~$(\sigma (\omega))_i= \omega_{i+1}$ for each $i\in \mathbb N$ 
  and $\omega= ( \omega_i)_{i\in \mathbb N} \in A^{\mathbb N}$). 
When a subset  $\Sigma$ of $A^{\mathbb N}$ is $\sigma$-invariant and  closed,  
we call it 
a \textit{subshift},
and $A$  the \textit{alphabet} of 
$\Sigma$.
When $\Sigma$ is of the form
\[
\Sigma =\{ (\omega_i)_{i\in\mathbb{N}}\in A^{\mathbb N} : \text{$M_{\omega_i \omega_{i+1}} =1 $ for all $i\in \mathbb N$}\}
\]
with  a matrix $M= (M_{ij})_{(i,j)\in A^2}$  each entry of which is $0$ or $1$, we call $\Sigma$ a \textit{Markov shift}.
When we emphasize the dependence of $\Sigma$ on $M$, it is denoted by $\Sigma _M$, and $M$ is called the \emph{adjacency matrix} of $\Sigma_M$.
In a similar manner, we define the shift map, a  subshift and  a Markov shift 
for the two-sided infinite product $A^{\mathbb Z}$ of $A$. 

For a subshift $\Sigma$ on an alphabet $A$, let 
$[u] := \left\{ (\omega_i)\in \Sigma : u=\omega_1\cdots\omega_n\right\}$ for each $u\in A^n$, $n\geq 1$,
and set
$
\mathcal L(\Sigma ) :=\left\{ u\in \bigcup _{n\geq 1} A^n : [u] \neq \emptyset \right\}.
$
We also denote $\mathcal{L}_n(\Sigma):=\{u\in\mathcal{L}(\Sigma):|u|=n\}$ for $n\ge 1$,
where $|u|$ means the length of $u$.
For $u=u_1\cdots u_m$ and $v=v_1\cdots v_n$ in $\mathcal{L}(\Sigma)$, we denote $uv=u_1\cdots u_mv_1\cdots v_n$.
Finally, we say that $\Sigma$ is \textit{transitive} if for any $u,v\in\mathcal{L}(\Sigma)$,
we can find $w\in\mathcal{L}(\Sigma)$ so that $uwv\in\mathcal{L}(\Sigma)$ holds.
In the rest of this paper, we denote by $h_{\rm top}(\Sigma)$ the \textit{topological entropy} of $\Sigma$,
and by $h(\mu)$ the \textit{metric entropy} of $\mu\in\mathcal{M}_{\sigma}(\Sigma)$.

\subsection{Markov diagram}

Let $T\colon X\to X$ be a piecewise monotonic map
and take an integer $k>1$ and $0=a_0<a_1<\cdots<a_k=1$ so that
$T|_{(a_{j-1},a_j)}$ is continuous and strictly monotone for each $1\le j\le k$.
Recall that $X_T:=\bigcap_{n=0}^{\infty}T^{-n}(\bigcup_{j=1}^k (a_{j-1},a_j))$.
We define the \textit{coding map} $\mathcal{I}\colon X_T\to \{1,\ldots,k\}^{\mathbb{N}}$ by
$$(\mathcal{I}(x))_i=j\text{ if and only if }T^{i-1}(x)\in (a_{j-1},a_j),$$
and denote the closure of $\mathcal{I}(X_T)$ by $\Sigma_T^+$.
Then it is well-known that $\Sigma_T^+$ is a subshift and $h_{\rm top}(\Sigma_T^+)=h_{\rm top}(T)$.
We call the system $(\Sigma_T^+,\sigma)$ a \textit{coding space} of $([0,1],T)$.

In what follows  we will construct Hofbauer's Markov Diagram, which is a countable oriented graph
with subsets of $\Sigma_T^+$ as vertices.
\begin{defi}
Let $C\subset\Sigma_T^+$ be a closed subset with $C\subset [j]$ for some $1\le j\le k$.
We say that a non-empty closed subset $D\subset\Sigma_T^+$ is a \textit{successor} of $C$ if
$D=[l]\cap\sigma(C)$ for some $1\le l\le k$.
\end{defi}

Now we define a set $\mathcal{D}_T$ of vertices by induction. First, we set
$\mathcal{D}_0:=\{[1],\ldots,[k]\}$. If $\mathcal{D}_n$ is defined for $n\ge 0$, then we set
$$\mathcal{D}_{n+1}:=\mathcal{D}_n\cup \{D:\text{there exists }C\in\mathcal{D}_n\text{ so that }D\text{ is a
successor of }C\}.$$
We note that $\mathcal{D}_n$ is a finite set for each $n\ge 0$ since the number of successors of any
closed subset of $\Sigma_T^+$ is
at most $k$ by the definition. Finally, we set
$$\mathcal{D}_T:=\bigcup_{n\ge 0}\mathcal{D}_n.$$
For the notational simplicity, we write $\mathcal{D}$ instead of $\mathcal{D}_T$
if no confusion arises.
The following lemma is obtained by an easy calculation and hence we omit the proof:
\begin{lem}
\label{ver}
(1) We have $\mathcal{D}=\{\sigma^{|u|-1}[u]:u\in\mathcal{L}(\Sigma_T^+)\}$.\\
(2) Let $u\in\mathcal{L}(\Sigma_T^+)$, $1\le j\le k$, and assume that $uj\in\mathcal{L}(\Sigma_T^+)$.
Then $\sigma^{|u|}[uj]$ is a successor of $\sigma^{|u|-1}[u]$.
\end{lem}
To get the oriented graph, which we call Hofbauer's Markov diagram, we insert an arrow from every
$C\in\mathcal{D}$ to all its successors. We write $C\rightarrow D$ to denote that
$D$ is a successor of $C$.
For $\mathcal{C}\subset \mathcal{D}$, we define a matrix
$M(\mathcal{C})=(M(\mathcal{C})_{C,D})_{(C,D)\in\mathcal{C}^2}$ by
$$M(\mathcal{C})_{C,D}=
\left\{
\begin{array}{ll}
1 & (\text{$C\rightarrow D$}), \\
0 & (\text{otherwise}).
\end{array}
\right.$$
Then $\Sigma_{M(\mathcal{C})}=\{(C_i)_{i\in\mathbb{Z}}\in\mathcal{C}^{\mathbb{Z}}:C_i\rightarrow
C_{i+1},i\in\mathbb{Z}\}$ is a Markov shift with a countable alphabet $\mathcal{C}$
and an adjacency matrix $M(\mathcal{C})$.
For the notational simplicity,
we denote $\Sigma_{\mathcal{C}}$ instead of $\Sigma_{M(\mathcal{C})}$ in the rest of this paper.

It is well-known that the topological entropy of $\Sigma_{\mathcal{D}}$ coincides with
that of $\Sigma_T^+$; $h_{\rm top}(\Sigma_{\mathcal{D}})=h_{\rm top}(\Sigma_T^+)$.
It is also known that there is a deep relationship between the Markov shift $\Sigma_{\mathcal{D}}$ and
a natural extension of $\Sigma_T^+$.
To be more precise, let
$$\Sigma_T:=\{(\omega_i)_{i\in\mathbb{Z}}\in\{1,\ldots,k\}^{\mathbb{Z}}:\omega_i\omega_{i+1}\cdots
\Sigma_T^+,i\in\mathbb{Z}\}$$
be a natural extension of $\Sigma_T^+$ and define a map
$\Psi\colon \Sigma_{\mathcal{D}}\to \{1,\ldots,k\}^{\mathbb{Z}}$ by
$$\Psi((C_i)_{i\in\mathbb{Z}}):=(\omega_i)_{i\in\mathbb{Z}}\text{ for }(C_i)_{i\in\mathbb{Z}}\in\Sigma_{\mathcal{D}},$$
where $1\le \omega_i\le k$ is the unique integer such that $C_i\subset [\omega_i]$ holds for each $i\in\mathbb{Z}$.
The following lemma, proved in \cite{H2}, states that $(\Sigma_{\mathcal{D}},\sigma)$ is topologically conjugate
to $(\Sigma_T,\sigma)$ except ``small" sets (see also \cite[Appendix]{Bu}):

\begin{lem}(\cite[Lemmas 2 and 3]{H2})
\label{hofbauer-extension}
The map $\Psi\colon\Sigma_{\mathcal{D}}\to\{1,\ldots,k\}^{\mathbb{Z}}$ is 
continuous and satisfies $\sigma\circ \Psi=\Psi\circ \sigma$. Moreover,
we can find two $\sigma$-invariant subsets $\overline{N}\subset\Sigma_{\mathcal{D}}$ and
$N\subset\Sigma_T$ satisfying the following properties:
\begin{itemize}
\item
We have
$\Psi(\Sigma_{\mathcal{D}}\setminus \overline{N})=\Sigma_T\setminus N$.

\item
The restriction map $\Psi\colon \Sigma_{\mathcal{D}}\setminus\overline{N}\to \Sigma_T\setminus N$
is bijective and bi-measurable.

\item
$N$ has no periodic point.

\item
For any invariant measure $\mu\in\mathcal{M}_{\sigma}(\Sigma_T)$ with $\mu(N)=1$, we have
$h(\mu)=0$.

\item
For any invariant measure $\overline{\mu}\in\mathcal{M}_{\sigma}(\Sigma_{\mathcal{D}})$ with
$\overline{\mu}(\overline{N})=1$, we have $h(\overline{\mu})=0$.
\end{itemize}
\end{lem}

Next, we prove that the condition (\ref{irreducible}) implies the transitivity of $\Sigma_{\mathcal{D}}$.
\begin{lem}
\label{transitive}
Suppose that $T$ satisfies the condition (\ref{irreducible}). Then
$\Sigma_{\mathcal{D}}$ is transitive.
\begin{proof}
Let $C,D\in\mathcal{D}$. To prove the lemma, it is sufficient to show that
there are $n\ge 1$ and $C_1,\ldots,C_n\in\mathcal{D}$ such that $C_1=C$, $C_n=D$ and
$C_i\rightarrow C_{i+1}$ for each $1\le i\le n-1$. By Lemma \ref{ver} (1), there are
$u=u_1\cdots u_l$ and $v=v_1\cdots v_m$ in $\mathcal{L}(\Sigma_T^+)$ such that
$C=\sigma^{l-1}[u]$ and $D=\sigma^{m-1}[v]$.
If we set $I:=\bigcap_{j=1}^lT^{-(j-1)}((a_{u_j-1},a_{u_j}))$, then
$I$ is a non-empty open interval since $u\in\mathcal{L}(\Sigma_T^+)$.
By the condition (\ref{irreducible}), there exists an open subinterval $L\subset I$,
$\tau\ge 1$, and $1\le i_0,\ldots,i_{\tau-1}\le k$ such that
$T^j(L)\subset (a_{i_j-1},a_{i_j})$ for $0\le j\le \tau-1$ and
$T^{\tau}(L)=(0,1)$. Since $T^{l-1}(I)\subset (a_{u_{l-1}},a_{u_l})$, we have $l\le \tau$.
We set
$$w=
\left\{
\begin{array}{ll}
i_l\cdots i_{\tau-1} & (\tau>l), \\
\lambda & (\tau=l),
\end{array}
\right.
$$
where $\lambda$ denotes the empty word.
Then it is not difficult to see that $\sigma^{\tau}[uw]=\Sigma_T^+$, which implies that
$\sigma^{\tau+m-1}[uwv]=D$.
Hence by Lemma \ref{ver} (2), we have the lemma.
\end{proof}
\end{lem}

In the rest of this section, we assume that $\Sigma_{\mathcal{D}}$ is transitive.
Combining \cite[Theorem 2 (iii)]{H2} and \cite[Page 377, Corollary 1 (ii)]{H3},
we have the following lemma:

\begin{lem}(\cite{H2,H3})
\label{mme}
There is a unique measure of maximal entropy both on $\Sigma _{\mathcal D}$ and $\Sigma _T^+$,
i.e., there is a unique ergodic measure $\overline{m}$ on $\Sigma _{\mathcal D}$ and
a unique ergodic measure $m^+$ on $\Sigma _T^+$ such that
$$h(\overline{m}) = h(m^+) = h_{\mathrm{top}} (\Sigma _T^+) . $$
Moreover, following properties hold for $\overline{m}$ and $m^+$.
\begin{enumerate}
\item
$m^+=\overline{m}\circ (\Psi^+)^{-1}$,
where $\Psi^+:=\pi\circ \Psi$ and $\pi\colon \Sigma_T\to \Sigma_T^+$ is
a natural projection; $(\omega_i)_{i\in\mathbb{Z}}\mapsto (\omega_i)_{i\in\mathbb{N}}$.

\item
There are two
families $\{L(C)\}_{C\in\mathcal{D}}$ and $\{R(C)\}_{C\in\mathcal{D}}$ of real positive numbers
satisfying the following properties:
\begin{itemize}
\item
$\overline{m}[C_1\cdots C_n]=L(C_1)R(C_n)\exp\{-nh_{\rm top}(\Sigma_T^+)\}$
for any $n\ge 1$ and $C_1\cdots C_n\in\mathcal{L}_n(\Sigma_{\mathcal{D}})$.

\item
Both $\sum_{C\in\mathcal{D}}L(C)$ and $\sup_{C\in\mathcal{D}}R(C)$ are finite.

\end{itemize}
\end{enumerate}
\end{lem}

Finally, we prove a weak Gibbs property for the unique measure of maximal entropy on $\Sigma_T^+$,
which plays a key role to obtain the lower large deviations bound in \S3.

\begin{lem}
\label{w-gibbs}
Let $m^+$ be the unique measure of maximal entropy on $\Sigma_T^+$
Then for any finite set $\mathcal{F}\subset\mathcal{D}$,
we can find $K=K_{\mathcal{F}}>1$ such that
\begin{equation}
\label{upper-gibbs}
m^+[u]\le K\exp\{-nh_{\rm top}(\Sigma_T^+)\}\ \ \ (n\ge 1,u\in\mathcal{L}_n(\Sigma_T^+)),
\end{equation}
$$m^+[u]\ge K^{-1}\exp\{-nh_{\rm top}(\Sigma_T^+)\}\ \ \ (n\ge 1,u\in\mathcal{L}_n(\Psi^+(\Sigma_{\mathcal{F}}))).$$
\begin{proof}
Let $\overline{m}$, $\{L(C)\}_{C\in\mathcal{D}}$ and $\{R(C)\}_{C\in\mathcal{D}}$ be as in Lemma \ref{mme}
For a finite set $\mathcal{F}\subset\mathcal{D}$, we put
$$K=K_{\mathcal{F}}:=\max\left\{2,LR,\left(\min_{C,D\in\mathcal{F}}L(C)R(D)\right)^{-1}\right\},$$
where we set $L:=\sum_{C\in\mathcal{D}}L(C)$ and $R:=\sup_{C\in\mathcal{D}}R(C)$.
Note that $K<\infty$ by Lemma \ref{mme}.
Let $n\ge 1$ and $u=u_1\cdots u_n\in\mathcal{L}_n(\Psi^+(\Sigma_{\mathcal{F}}))$. Then we can find a word
$C_1\cdots C_n\in\mathcal{L}_n(\Sigma_{\mathcal{F}})$ such that
$ [C_1\cdots C_n]\subset(\Psi^+)^{-1}[u]$.
Hence it follows from Lemma \ref{mme} that
$$m^+[u]\ge \overline{m}[C_1\cdots C_n]=L(C_1)R(C_n)\exp\{-nh_{\rm top}(\Sigma_T^+)\}
\ge K^{-1}\exp\{-nh_{\rm top}(\Sigma_T^+)\}.$$

We show (\ref{upper-gibbs}). First, we note that $\Psi^+(\Sigma_{\mathcal{D}})=\Sigma_T^+$
(see \cite{H3} for instance).
For $n\ge 1$ and $u\in\mathcal{L}_n(\Sigma_T^+)$,
we set
$$\mathcal{C}=\mathcal{C}(u):=\{C\in\mathcal{D}:[C]\cap(\Psi^+)^{-1}[u]\not=\emptyset\}.$$
Then for any $C\in\mathcal{C}$, we can find a unique word
$P(C)=C_1\cdots C_n\in\mathcal{L}_n(\Sigma_{\mathcal{D}})$
such that $C_1=C$ and $[C]\cap (\Psi^+)^{-1}[u]=[P(C)]$.
Then it is easy to see that $\bigcup_{C\in\mathcal{D}}[P(C)]=(\Psi^+)^{-1}[u]$.
Therefore we have
\begin{eqnarray*}
m[u]
&\le&
\sum_{C\in\mathcal{C}}\overline{m}[P(C)]\\
&\le&
\sum_{C\in\mathcal{D}}L(C)R\exp\{-nh_{\rm top}(\Sigma_T^+)\}\\
&\le&
K\exp\{-nh_{\rm top}(\Sigma_T^+)\},
\end{eqnarray*}
which proves the lemma.
\end{proof}
\end{lem}

\section{Proofs of Theorems \ref{mod1} and \ref{interval}}
In this section, we give proofs of Theorems \ref{mod1} and \ref{interval}.
First, we show the following analogous result of Theorem \ref{interval} on coding spaces:
\begin{maintheorem}
\label{mainprop}
Let $T\colon X\to X$ be a transitive piecewise monotonic map, $(\Sigma_T^+,\sigma)$ be a coding space of
$(X,T)$, and
$m^+$ be the unique measure of maximal entropy on $\Sigma_T^+$.
Suppose that $\mathcal{M}_{\sigma}^p(\Sigma_T^+)$ is dense in $\mathcal{M}^e_{\sigma}(\Sigma_T^+)$
and $\Sigma_{\mathcal{D}_T}$ is transitive.
Then $(\Sigma_T^+,\sigma)$ satisfies the large deviation principle with $m^+$ and the rate function
$\mathcal{J}^+\colon\mathcal{M}(\Sigma_T^+)\to [0,\infty]$ is given by
\begin{equation}
\label{code-rate}
\mathcal{J}^+(\mu^+)=
\left\{
\begin{array}{ll}
h_{\rm top}(\Sigma_T^+)-h(\mu^+) & (\mu^+\in \mathcal{M}_{\sigma}(\Sigma_T^+)), \\
\infty & (\text{otherwise}).
\end{array}
\right.
\end{equation}
\begin{proof}
To prove the theorem, it is sufficient to show the following:
\begin{equation}
\label{upper}
\limsup_{n\rightarrow\infty}\frac{1}{n}\log m^+((\delta_n^{\sigma})^{-1}(\mathcal{K}))\le
-\inf_{\mu^+\in\mathcal{K}\cap\mathcal{M}_{\sigma}(\Sigma_T^+)}(h_{\rm top}(\Sigma_T^+)-h(\mu^+))
\end{equation}
holds for any closed set $\mathcal{K}\subset\mathcal{M}(\Sigma_T^+)$ and
\begin{equation}
\label{lower}
\liminf_{n\rightarrow\infty}\frac{1}{n}\log m^+((\delta_n^{\sigma})^{-1}(\mathcal{U}))\ge
-\inf_{\mu^+\in\mathcal{U}\cap\mathcal{M}_{\sigma}(\Sigma_T^+)}(h_{\rm top}(\Sigma_T^+)-h(\mu^+))
\end{equation}
holds for any open set
$\mathcal{U}\subset\mathcal{M}(\Sigma_T^+)$.
We note that
the equation (\ref{upper-gibbs}) in Lemma \ref{w-gibbs} enables us to show (\ref{upper})
in a similar way to the proof of \cite[Theorem A, \S4]{CTY}.

In what follows we will show (\ref{lower}).
Let $\epsilon>0$, $\mu^+\in\mathcal{M}_{\sigma}(\Sigma_T^+)$ and $\mathcal{U}$ be an open neighborhood
of $\mu^+$ in $\mathcal{M}(\Sigma_T^+)$.
To show (\ref{lower}), it is enough to show that
$$\liminf_{n\rightarrow\infty}\frac{1}{n}\log m^+((\delta_n^{\sigma})^{-1}(\mathcal{U}))\ge
h(\mu^+)-h_{\rm top}(\Sigma_T^+)-2\epsilon.$$
First, we prove the following proposition:
\begin{prop}
\label{e-dense}
There exist a finite set $\mathcal{F}\subset \mathcal{D}$ and
$\rho^+\in\mathcal{M}_{\sigma}^e(\Psi^+(\Sigma_{\mathcal{F}}))$ such that
$\rho^+\in \mathcal{U}$ and $h(\rho^+)\ge h(\mu^+)-\epsilon$.
\begin{proof}
If $h(\mu^+)-\epsilon\le 0$, then we choose $\rho^+\in\mathcal{M}_{\sigma}^p(\Sigma_T^+)\cap \mathcal{U}$
and take a periodic point $\omega\in\Sigma_T^+$ in the support of $\rho^+$.
It follows from \cite[Theorem 8]{H3} that there are finite vertices $C_1,\ldots, C_n\in\mathcal{D}$ so that
$\Psi^+(\cdots C_1\cdots C_nC_1\cdots C_n\cdots)=\omega$ holds.
We set $\mathcal{F}:=\{C_1,\ldots,C_n\}$. Then it is clear that
$\rho^+\in\mathcal{M}_{\sigma}^e(\Psi^+(\Sigma_{\mathcal{F}}))$ and $h(\rho^+)=0\ge h(\mu^+)-\epsilon$.
This suffices for the case $h(\mu^+)-\epsilon\le 0$.

Hereafter assume that $h(\mu^+)-\epsilon>0$.
First, we define two push forward maps $\pi_{\ast}\colon \mathcal{M}_{\sigma}(\Sigma_T)
\to\mathcal{M}_{\sigma}(\Sigma_T^+)$ and
$\Psi_{\ast}\colon\mathcal{M}(\Sigma_{\mathcal{D}})\to
\mathcal{M}(\Sigma_T)$
by $\pi_{\ast}(\vartheta):=\vartheta\circ\pi^{-1}$ for $\vartheta\in\mathcal{M}_{\sigma}(\Sigma_T)$, and
$\Psi_{\ast}(\overline{\vartheta}):=\overline{\vartheta}\circ\Psi^{-1}$ for
$\overline{\vartheta}\in\mathcal{M}(\Sigma_{\mathcal{D}})$.
Then it is well-known that $\pi_{\ast}$ is a homeomorphism and
$h(\vartheta)=h(\pi_{\ast}(\vartheta))$ for any $\vartheta\in\mathcal{M}_{\sigma}(\Sigma_T)$.
Moreover, it follows from Lemma \ref{hofbauer-extension} that $\Psi_{\ast}$ is continuous and
the restriction map $\Psi_{\ast}\colon\mathcal{M}_{\sigma}(\Sigma_{\mathcal{D}}\setminus\overline{N})\to
\mathcal{M}_{\sigma}(\Sigma_T\setminus N)$
is well-defined, bi-measurable, and
$h(\overline{\vartheta})=h(\Psi_{\ast}(\overline{\vartheta}))$ holds for any $\overline{\vartheta}\in
\mathcal{M}_{\sigma}(\Sigma_{\mathcal{D}}\setminus\overline{N})$.
Here $N$ and $\overline{N}$ be as in Lemma \ref{hofbauer-extension}.

We set $\mu:=\pi_{\ast}^{-1}(\mu^+)\in\pi_{\ast}^{-1}(\mathcal{U})$.
Then we can find $0\le c\le 1$ and $\mu_1,\mu_2\in\mathcal{M}_{\sigma}(\Sigma_T)$
such that
$\mu=c\mu_1+(1-c)\mu_2$, $\mu_1(\Sigma_T\setminus N)=1$ and $\mu_2(N)=1$ hold.
We note that $h(\mu_2)=0$ by Lemma \ref{hofbauer-extension}.
Since $\mathcal{M}^p_{\sigma}(\Sigma_T^+)$ is dense in $\mathcal{M}_{\sigma}^e(\Sigma_T^+)$
by the assumption, it follows from \cite[Theorem A]{Y2} that
$\mathcal{M}_{\sigma}^p(\Sigma_T^+)$ is also dense in $\mathcal{M}_{\sigma}(\Sigma_T^+)$.
By the definition of $\pi_{\ast}$, one can easily see that
$\mathcal{M}^p_{\sigma}(\Sigma_T^+)=\pi_{\ast}(\mathcal{M}^p_{\sigma}(\Sigma_T))$
and hence $\mathcal{M}^p_{\sigma}(\Sigma_T)$ is also dense in $\mathcal{M}_{\sigma}(\Sigma_T)$.
Therefore, we can find $\nu_2\in\mathcal{M}_{\sigma}^p(\Sigma_T)$
such that $\nu:=c\mu_1+(1-c)\nu_2\in\pi_{\ast}^{-1}(\mathcal{U})$
and $h(\nu)=h(\mu)$.
Since $N$ has no periodic point, $\nu_2(\Sigma_T\setminus N)=1$, which also implies that
$\nu(\Sigma_T\setminus N)=1$.
Hence, we can define $\overline{\nu}:=\Psi_{\ast}^{-1}(\nu)\in (\pi_{\ast}\circ\Psi_{\ast})^{-1}(\mathcal{U})$.
Note that $\overline{\nu}$ is supported on a transitive countable Markov shift $\Sigma_{\mathcal{D}}$ and
$(\pi_{\ast}\circ\Psi_{\ast})^{-1}(\mathcal{U})$ is open.
Therefore, it follows from \cite[Main Theorem]{T} that there exist
a finite set $\mathcal{F}\subset \mathcal{D}$ and an ergodic measure
$\overline{\rho}\in (\pi_{\ast}\circ\Psi_{\ast})^{-1}(\mathcal{U})$ supported on $\Sigma_{\mathcal{F}}$
such that $h(\overline{\rho})\ge h(\overline{\nu})-\epsilon>0$.
Since $h(\overline{\rho})>0$, we have
$\overline{\rho}(\overline{N})=0$ by Lemma \ref{hofbauer-extension}.
We define an ergodic measure $\rho^+$ on $\Sigma_T^+$ by $\rho^+:=\pi_{\ast}(\Phi_{\ast}(\overline{\rho}))$.
Clearly, we have $\rho^+\in \mathcal{M}_{\sigma}^e(\Psi^+(\Sigma_{\mathcal{F}}))\cap \mathcal{U}$.
Moreover, since $\overline{\rho}\in\mathcal{M}_{\sigma}(\mathcal{M}(\Sigma_{\mathcal{D}}\setminus\overline{N}))$,
we have $h(\overline{\rho})=h(\rho^+)$, which implies that $h(\rho^+)\ge h(\mu^+)-\epsilon$.
\end{proof}
\end{prop}

We continue the proof of Theorem \ref{mainprop}.
Note that $\Psi^+(\Sigma_{\mathcal{F}})\subset\Sigma_T^+$ is a subshift and
$\rho^+(\Psi^+(\Sigma_{\mathcal{F}}))=1$.
Hence combining \cite[Propositions 2.1 and 4.1]{PS}, we have the following:

\begin{lem}(\cite{PS})
There exists $M\in\mathbb{N}$ such that for any $n\ge M$,
$$\#\{u\in\mathcal{L}_n(\Psi^+(\Sigma_{\mathcal{F}})):[u]\subset (\delta_n^{\sigma})^{-1}(\mathcal{U})\}\ge
\exp\{n(h(\rho^+)-\epsilon)\}.$$
\end{lem}

For the notational simplicity, for each $n\ge M$, we set
$\mathcal{L}_{n,\rho^+}:=\{u\in\mathcal{L}_n(\Psi^+(\Sigma_{\mathcal{F}})):[u]\subset (\delta_n^{\sigma})^{-1}(\mathcal{U})\}.$
Since $\mathcal{F}$ is finite, it follows from Lemma \ref{w-gibbs} that there is $K>1$ such that
$m^+[u]\ge K^{-1}\exp\{-nh_{\rm top}(\Sigma_T^+)\}$
for any $u\in\mathcal{L}_{n,\rho^+}$ and $n\ge M$.
Therefore, by $\bigcup_{u\in\mathcal{L}_{n,\rho^+}}[u]\subset (\delta_n^{\sigma})^{-1}(\mathcal{U})$,
we have
\begin{eqnarray*}
{}
& &
\liminf_{n\rightarrow\infty}\frac{1}{n}\log m^+((\delta_n^{\sigma})^{-1}(\mathcal{U}))\\
&\ge&
\liminf_{n\rightarrow\infty}\frac{1}{n}\log m^+\left(\bigcup_{u\in\mathcal{L}_{n,\rho^+}}[u]\right) \\
&=&
\liminf_{n\rightarrow\infty}\frac{1}{n}\log \sum_{u\in\mathcal{L}_{n,\rho^+}}m^+[u]\\
&\ge&
\liminf_{n\rightarrow\infty}\frac{1}{n}\log \#\mathcal{L}_{n,\rho^+}K^{-1}
\exp\{-nh_{\rm top}(\Sigma_T^+)\}\\
&\ge&
\liminf_{n\rightarrow\infty}\frac{1}{n}\log K^{-1}\exp\{n(h(\rho^+)-h_{\rm top}(\Sigma_T^+)-\epsilon)\}\\
&=&
h(\mu^+)-h_{\rm top}(\Sigma_T^+)-2\epsilon,
\end{eqnarray*}
which proves Theorem \ref{mainprop}.

\end{proof}
\end{maintheorem}

\noindent
{\bf Proof of Theorem \ref{interval}.}
Let $T\colon X\to X$ and $m$ be as in Theorem \ref{interval} and
assume that $\mathcal{M}_T^p(X)$ is dense in $\mathcal{M}_T^e(X)$.
Then it follows from \cite[Theorem A]{Y2} that $\mathcal{M}_{\sigma}^p(\Sigma_T^+)$
is dense in $\mathcal{M}_{\sigma}(\Sigma_T^+)$.
Since $T$ satisfies the condition (\ref{irreducible}), $\Sigma_{\mathcal{D}}$ is transitive by
Lemma \ref{transitive}.
Hence it follows from Theorem \ref{mainprop}
that $(\Sigma_T^+,\sigma)$ satisfies the large deviation principle with
the unique measure of maximal entropy $m^+$ on $\Sigma_T^+$ and the rate function
$\mathcal{J}^+\colon\mathcal{M}(\Sigma_T^+)\to [0,\infty]$ is given by the equation (\ref{code-rate}).

Define a map $\Phi\colon\Sigma_T^+\to X$ by
$\Phi(\omega):=y$, where $y$ is a unique element of a set
$\bigcap_{n\ge 0}{\rm cl}(T^{-n}(a_{\omega_{n+1}-1},a_{\omega_{n+1}}))$
and let $\Phi_{\ast}\colon\mathcal{M}(\Sigma_T^+)\to \mathcal{M}([0,1])$
be a push forward map; $\mu^+\mapsto \mu^+\circ\Phi^{-1}$.
Then it is well-known that $\Phi$ is a continuous surjection,
$\Phi\circ\sigma(\omega)=T\circ\Phi(\omega)$ for $\omega\in\Phi^{-1}(X_T)$, $m=\Phi_{\ast}(m^+)$,
and $h(\mu^+)=h_{\Phi_{\ast}(\mu^+)}(T)$ if $\Phi_{\ast}(\mu^+)\in\mathcal{M}_T(X)$
(see \cite{F} for instance).
Therefore, it follows from Contraction Principle (\cite[Theorem 4.2.1]{DZ}) that
$(X,T)$ satisfies the large deviation principle with $m$ and the rate function
$\mathcal{I}\colon\mathcal{M}(X)\to [0,\infty]$ is given by
$$\mathcal{J}(\mu):=\inf\{\mathcal{J}^+(\mu^+):\mu^+\in\mathcal{M}(\Sigma_T^+),\mu=\Phi_{\ast}(\mu^+)\}.$$
Since $\mathcal{M}_{\sigma}^p(\Sigma_T^+)$ is dense in $\mathcal{M}_{\sigma}(\Sigma_T^+)$,
it follows from \cite[Lemma 2.6]{Y2} that $\Phi_{\ast}(\mathcal{M}_{\sigma}(\Sigma_T^+))={\rm cl}(\mathcal{M}_T(X_T))$.
Hence $\mathcal{J}(\mu)=\infty$ holds on $\mathcal{M}(X)\setminus {\rm cl}(\mathcal{M}_T(X_T))$
and the set
$\{h_{\rm top}(\Sigma_T^+)-h(\mu^+):\mu^+\in\mathcal{M}_{\sigma}(\Sigma_T^+),\mu=\Phi_{\ast}(\mu^+)\}$
is not empty for any $\mu\in {\rm cl}(\mathcal{M}(X_T))$. Therefore, for any $\mu\in {\rm cl}(\mathcal{M}(X_T))$,
we have
\begin{eqnarray*}
\mathcal{J}(\mu)
&=&
\inf\{\mathcal{J}^+(\mu^+):\mu^+\in\mathcal{M}_{\sigma}(\Sigma_T^+),\mu=\Phi_{\ast}(\mu^+)\} \\
&=&
\inf\{h_{\rm top}(\Sigma_T^+)-h(\mu^+):\mu^+\in\mathcal{M}_{\sigma}(\Sigma_T^+),\mu=\Phi_{\ast}(\mu^+)\} \\
&=&
h_{\rm top}(T)-h_{\mu}(T),
\end{eqnarray*}
which proves Theorem \ref{interval}.
\vspace{0.3cm}\\
{\bf Proof of Theorem \ref{mod1}.}
Let $0\le \alpha <1$, $\beta>2$, and $T=T_{\alpha,\beta}$ be as in Theorem \ref{mod1}.
Then it follows from \cite[Theorem 2]{H} that
$\mathcal{M}_{\sigma}^p(\Sigma_T^+)$ is dense in $\mathcal{M}_{\sigma}(\Sigma_T^+)$,
and hence $\mathcal{M}_T^p(X)$ is also dense in $\mathcal{M}_T^e(X)$.
It follows from \cite[Proposition 3.14]{CLR} that $T$ satisfies
the condition (\ref{irreducible}).
Hence by Theorem \ref{interval}, $(X,T)$ satisfies
the large deviation principle with the unique measure of maximal entropy $m$,
and the rate function $\mathcal{J}\colon\mathcal{M}(X)\to [0,\infty]$ is given by (\ref{rate2}).
Finally, since $T$ is piecewise increasing and right continuous,
it follows from \cite[Theorem B and Lemma 2.7]{Y2} that ${\rm cl}(\mathcal{M}_T(X_T))=\mathcal{M}_T(X)$,
which proves Theorem \ref{mod1}.

\section{Further applications}


In this section, we apply Theorem \ref{interval} for $(-\beta)$-transformations with $(1+\sqrt{5})/2<\beta<2$,
and give an example of piecewise monotonic maps $T$ to which Theorem \ref{interval} is applicable,
and $\mathcal{M}_T(X)$ does not coincide  with ${\rm cl}(\mathcal{M}_T(X_T))$.

\begin{example}
\label{counter}
Let $(1+\sqrt{5})/2<\beta<2$ and define a \textit{$(-\beta)$-transformation}
$T=T_{-\beta}\colon X\to X$ by
$$T_{-\beta}(x)=
\left\{
\begin{array}{ll}
-\beta x+1 & (0\le x<\frac{1}{\beta}) \\
-\beta x+2 & (\frac{1}{\beta}\le x\le 1).
\end{array}
\right.$$
Since $\beta>(1+\sqrt{5})/2$, we can prove that $T$ satisfies the condition (\ref{irreducible})
in a similar way to the proof of \cite[Proposition 8]{G}, where the exactness is shown for the absolutely
continuous $T$-invariant measure.
On the other hand, it is clear that $T$ has two monotonic pieces, i.e.,
there is $0<a<1$ such that both $T|_{(0,a)}$ and $T|_{(a,1)}$ are continuous and strinctly monotone.
Since $T$ is transitive,  it follows from \cite[Theorem 2]{HR} that $\mathcal{M}^p_{\sigma}(\Sigma_T^+)$
is dense in $\mathcal{M}_{\sigma}(\Sigma_T^+)$, and hence
$\mathcal{M}_T^p(X)$ is also dense in $\mathcal{M}^e_T(X)$.
Therefore, by Theorem \ref{interval},
$(X,T)$ satisfies the large deviation principle with the unique measure of maximal entropy $m$,
and the rate function $\mathcal{J}\colon \mathcal{M}(X)\to [0,\infty]$ is given by (\ref{rate2}).

Hereafter, let $\beta=1.7548\cdots$ be a unique real root of the equation
$t^3-2t^2+t-1=0$.
Then it is proved in \cite[Example 3.1]{Y2} that $\mathcal{M}_T^p(X)$ does not dense in $\mathcal{M}_T(X)$.
Hence by \cite[Lemmas 2.6 and 2.7]{Y2},
we have ${\rm cl}(\mathcal{M}_T(X_T))\not=\mathcal{M}_T(X)$.
\end{example}

\vspace*{3mm}

\noindent
\textbf{Acknowledgement.}~ 
The authors would like to express their gratitude to Hiroki Takahasi for fruitful discussions.
The first author was partially supported by JSPS KAKENHI Grant Number 16K05179 and
the second author was partially supported by JSPS KAKENHI Grant Number 18K03359.

\end{document}